\documentclass{article}
\usepackage{latexsym,amssymb}
\begin{document}
\begin{center}
{\bf On the minimal norm of a non-regular generalized character of an arbitrary finite group}\\
\emph{Geoffrey R. Robinson\\
Department of Mathematics\\
University of Aberdeen\\
Aberdeen AB24 3UE\\
Scotland\\
United Kingdom\\
e-mail:g.r.robinson@abdn.ac.uk}
\end{center}

\medskip
\noindent {\bf INTRODUCTION AND NOTATION:} As usual, for $G$ a
finite group, we denote the set of non-identity elements of $G$ by
$G^{\#},$ the number of conjugacy classes of $G$ by $k(G)$ and the
set of complex irreducible characters of $G$ by ${\rm Irr}(G).$ When
$H$ is a non-trivial subgroup of $G,$ we let $k_{+}(G,H)$ denote the
number of irreducible characters of $G$ which do not vanish
identically on $H^{\#}$\\ (and we remark that $k_{+}(G,H) = k(G)$
when $|H| \geq \sqrt{|G|}$).

\medskip
For $\theta$ a generalized character of a finite group $G,$ we set
$$\sigma(G^{\#}, \theta) = \sum_{g \in G^{\#}} |\theta(g)|^{2}.$$
Notice that $$\sigma(G^{\#}, \theta) = |G|\langle \theta,\theta
\rangle - \theta(1)^{2},$$ where $\langle ,\rangle$ is the usual
inner product on the space of complex-valued class functions of $G.$
We denote the regular character of $G$ by $\rho.$ We let ${\rm
cf}(G)$ denote the space of complex-valued class functions of $G,$
and for $\alpha,\beta \in {\rm cf(G)},$ we set
$$(\alpha,\beta) = |G|\langle \alpha,\beta \rangle -
\alpha(1)\beta(1).$$ Then $(,)$ induces an inner product on the
quotient space ${\rm cf}(G)/\mathbb{C}\rho.$ We label the
irreducible characters of $G$ as $\mu_{1} ( = 1),\mu_{2},\ldots,
\mu_{k}$ where $\mu_{1}(1) \leq \ldots \leq \mu_{k}(1).$

\medskip
In the earlier papers [1],[2] and [3], we considered the minimum
possible value of $\sigma(G^{\#}, \theta)$ for $\theta$ a
generalized character (other than a multiple of the regular
character) of a non-trivial finite group $G$. We denote this minimum
possible value by $m(G)$ (which is clearly a positive integer). Note
that with the above notation, this is the minimum non-zero value of
$(\theta,\theta)$ for $\theta \in \mathbb{Z}{\rm Irr}(G) \backslash
\mathbb{Z}\rho.$

\medskip
We proved in [1] that $m(G) \geq k(G)-1$ for every non-trivial
finite group $G.$ In the papers [2],[3] we proved that $m(G) = |G| -
\mu_{k}(1)^{2}$ whenever $G$ is a non-trivial nilpotent group. In
particular, in that case,
$$m(G) \geq |G| - [G:Z(G)]$$ whenever $G$ since $\mu_{k}(1)^{2}$
divides $[G:Z(G)].$

\medskip
We asked in [3] whether it was the case that for every non-trivial
finite group $G$ we have $$m(G) \geq \frac{|G|}{\mu_{k}(1)} - 1 $$
(which would yield $m(G) > \sqrt{|G|}-1)$ . We answer that question
in the affirmative here. Notice that if $G$ is a Frobenius group
with Abelian kernel $K$ and complement $H,$ then $m(G) = |K| - 1.$
In that case, the maximal degree of an irreducible character of $G$
is $|H|,$ and we do have
$$m(G) = \frac{|G|}{\mu_{k}(1)} - 1.$$

\medskip
For if $\theta$ is a generalized character of $G$ which does not
vanish identically on $H^{\#},$ we see that $\sigma(G^{\#}, \theta)
\geq |K|\sigma(H^{\#}, {\rm Res}^{G}_{H}(\theta)) \geq |K|,$ as
there are $|K|$ conjugates of $H$ which intersect trivially with
each other. On the other hand, if the generalized character $\theta$
does not vanish identically on $K^{\#},$ then $$\sigma(G^{\#},
\theta) \geq m(K) = |K|-1,$$ as $K$ is Abelian. Finally, let
$\theta$ be the sum of the distinct irreducible characters of $G$
which do not contain $K$ in their kernels. By the orthogonality
relations, we see that $\theta(x) = -1$ for all $x \in K^{\#},$
while also $\theta$ vanishes identically outside $K,$ so that
$\sigma(G^{\#},\theta) = |K|-1,$ as claimed.

\medskip
In view of these examples, our main theorem may be regarded as
optimal in this generality. We prove:

\medskip
\noindent {\bf THEOREM 1 :} \emph{ Let $G$ be a non-trivial finite
group. Then
$$m(G) \geq \frac{|G|}{\mu_{k}(1)} - 1.$$ Furthermore, equality is
attained if and only if $G$ is either Abelian, or else is a
Frobenius group with an Abelian Frobenius kernel. }

\medskip
As noted above, Theorem 1 has an interpretation in terms of the
minimal non-zero value of a certain semi-definite integral quadratic
form.

\medskip
\noindent {\bf COROLLARY 2:} \emph{ Let $G$ be a non-trivial finite
group and let $a_{1},a_{2},\ldots, a_{k}$ be any $k$ integers such
that for some $i \neq j,$ we have\\ $a_{i}\mu_{j}(1) \neq
a_{j}\mu_{i}(1).$ Then
$$\sum_{i < j}(a_{i}\mu_{j}(1) -a_{j}\mu_{i}(1))^{2} \geq
\frac{|G|}{\mu_{k}(1)} - 1.$$}

\medskip
\noindent {\bf COROLLARY 3:} \emph{ Let $H$ be a non-trivial
subgroup of the finite group $G$ and let $\mu$ be a complex
irreducible character of maximal degree of $H.$ Then we have
$$k_{+}(G,H) \leq \frac{\mu(1)|H|-\mu(1)}{|H|-\mu(1)}{\rm max}_{h
\in H^{\#}} |C_{G}(h)|.$$}

\begin{center}
{\bf PROOF OF THEOREM 1}
\end{center}
Suppose that $m(G) < \frac{|G|}{\mu_{k}(1)}.$ To prove the Theorem,
by the cases discussed above, we may suppose that $G$ is
non-Abelian. Let us choose a generalized character $\theta$ of $G$
such that $m(G) = \sigma(G^{\#},\theta).$ We will prove that
$\theta(x) \in \{0,-1\}$ for all $x \in G^{\#},$ and that we may
suppose that $\theta$ is a multiplicity free irreducible character
of $G$ whose irreducible constituents all have degree $\mu_{k}(1),$
and that if there are $r$ such constituents, then $|G| = \mu_{k}(1)
+ r\mu_{k}(1)^{2}.$

\medskip
If $\theta$ has the form $t\rho + \delta \mu_{i}$ for some integer
$t,$ some sign $\delta$ and some $i,$ then we may replace $\theta$
by $\mu_{i}.$ Then we have
$$\sigma(G^{\#},\theta) = |G| - \mu_{i}(1)^{2},$$ so we may suppose that $i = k,$
as $\mu_{i}$ must also have maximal degree by the choice of
$\theta.$ We then have
$$\sigma(G^{\#},\theta) = |G| - \mu_{k}(1)^{2} \geq
\frac{|G|}{\mu_{k}(1)} - 1.$$ Furthermore, equality forces either
$\mu_{k}(1) = 1$ (in which case $G$ is Abelian, contrary to
hypothesis), or else $|G| = \mu_{k}(1)^{2} + \mu_{k}(1).$ In the
latter case, we have $\sum_{i=1}^{k-1}\mu_{i}(1)^{2} = \mu_{k}(1).$
Now $\langle \theta,\mu_{1} \rangle = 0,$ so that $\sum_{g \in
G^{\#}} \mu_{k}(g) = -\mu_{k}(1).$ Under the current assumptions, we
also have $$\sigma(G^{\#},\theta) = \mu_{k}(1).$$ Now $\mu_{k}$ is
the unique irreducible character of $G$ of its degree, so it is
rational-valued. Hence we must have $\mu_{k}(g) = -1$ whenever
$\mu_{k}(G) \neq 0.$ In order to establish the initial claim, then,
we may now suppose that $\theta$ does not have the form $t\rho +
\delta \mu_{i}$ for some integer $t,$ some sign $\delta$ and some
$i.$

\medskip
By the arguments of  [2] and [3], we may suppose that $\theta$ is a
rational-valued character of $G$ with $0 < \theta(1) <
\frac{|G|}{2}$ ( in which irreducible characters of the same degree
occur with equal multiplicity), and that $\langle \theta,\mu_{i}
\rangle$ is the closest integer to $\frac{\theta(1)\mu(1)}{|G|}$ for
each irreducible character $\mu_{i}$ of $G.$ For the convenience of
the reader, we review the necessary results using the notation
currently adopted. After subtracting integer multiples of $\rho$
from $\theta$ and replacing $\theta$ by $-\theta$ if necessary, we
may suppose that $0 < \theta(1) \leq  \frac{|G|}{2}.$

\medskip
For (by the choice of $\theta$) we have $$(\theta - \delta
\mu_{i},\theta- \delta \mu_{i}) \geq (\theta,\theta)$$ for each $i$
and any choice of sign $\delta.$ This leads easily to
$$|(\theta,\mu_{i})| \leq \frac{|G| - \mu_{i}(1)^{2}}{2}$$
and then $$\left| \langle \theta, \mu_{i} \rangle -
\frac{\theta(1)\mu_{i}(1)}{|G|} \right| \leq \frac{1}{2} \left( 1-
\frac{\mu_{i}(1)^{2}}{|G|} \right).$$ Hence $\langle \theta,
\mu_{i}\rangle$ is the (unique under current assumptions) closest
integer to $\frac{\theta(1)\mu_{i}(1)}{|G|},$ and only depends on
$\mu_{i}(1).$ Notice that we can't now have $\theta(1) =
\frac{|G|}{2},$ for if we did, we would have the contradiction
$$\left| \langle \theta, \mu_{1} \rangle - \frac{1}{2} \right| <
\frac{1}{2}.$$ Thus we have $0 < \theta(1) < \frac{|G|}{2}.$ Also,
$\theta$ is a rational (integer) valued character, since it contains
algebraically conjugate irreducible characters with the same
non-negative multiplicity.

\medskip
Notice that the trivial character does not now occur as an
irreducible constituent of $\theta.$ Hence we have $\sum_{g \in
G^{\#}} \theta(g) = -\theta(1),$ so that
$$\sum_{g \in G^{\#}} |\theta(g)| \geq \theta(1)$$ and certainly
$$\sigma(G^{\#},\theta) \geq \theta(1)$$ ( as $\theta$ is
integer-valued).

\medskip
Now let us write $$\theta = \sum_{i=1}^{k} a_{i}\mu_{i},$$ and set
$$\epsilon_{i} = a_{i} - \frac{\theta(1)\mu_{i}(1)}{|G|}.$$ Then we
have $$\sum_{i=1}^{k} \epsilon_{i}\mu_{i}(1) = 0$$ and
$$\sigma(G^{\#},\theta) = |G| \left(\sum_{i=1}^{k}
|\epsilon_{i}|^{2}\right)$$ ( for notice that $\theta$ agrees with
$\sum_{i=1}^{k} \epsilon_{i}\mu_{i}$ on $G^{\#}).$ We have noted
already that $|\epsilon_{i}| < \frac{1}{2}$ for each $i.$

\medskip
Suppose that $a_{i} \geq 2$ for some $i.$ Then we have
$$\frac{\theta(1)\mu_{i}(1)}{|G|} > \frac{3}{2},$$ so that $$\sigma(G^{\#},\theta)
\geq \theta(1)
> \frac{G}{\mu_{i}(1)} \geq \frac{|G|}{\mu_{k}(1)},$$ contrary to
hypothesis. Hence we have $a_{i} \in \{0,1\}$ for each $i.$

\medskip
Now suppose that $\epsilon_{i} < 0$ for $i = 1,\ldots,k-r,$ and that
$\epsilon_{i} > 0$ for $i = k-r+1,\ldots,k.$ Then
$$\theta = \sum_{i=1}^{r} \mu_{k-r + i}.$$ Since we have dealt with
the case that $\theta$ is irreducible, we may, and do, suppose that
$r > 1.$

\medskip
By assumption, we have $$\frac{|G|}{\mu_{k}(1)} - 1 \geq m(G) =
\sigma(G^{\#},\theta) \geq \theta(1),$$ while we clearly also have
$$\theta(1) \geq \frac{|G| - \sum_{i=1}^{k-r}
\mu_{i}(1)^{2}}{\mu_{k}(1)},$$ so that $$\sum_{i=1}^{k-r}
\mu_{i}(1)^{2} \geq  \mu_{k}(1)$$ ( and hence $$|G| \geq \mu_{k}(1)
+ \sum_{i=k-r+1}^{k} \mu_{i}(1)^{2}).$$

\medskip
Now we have
$$\sigma(G^{\#},\theta) = \sum_{i < j} (a_{i}\mu_{j}(1) -
a_{j}\mu_{i}(1))^{2} \geq r \left(\sum_{i=1}^{k-r}
\mu_{i}(1)^{2}\right),$$ using the fact that $a_{j} = 1$ for $j >
k-r$ and $a_{i} = 0$ for $i \leq k-r.$

\medskip
Thus we have (by assumption) $$r\left[|G|- \left(\sum_{i=k-r+1}^{k}
\mu_{i}(1)^{2}\right)\right] \leq \sigma(G^{\#},\theta) \leq
\frac{|G|}{\mu_{k}(1)} -1,$$ so that
$$|G|\left( 1 - \frac{1}{r\mu_{k}(1)}\right) <\left(\sum_{i=k-r+1}^{k}
\mu_{i}(1)^{2}\right).$$

\medskip
Now for $n \geq 2,$ we have $$\frac{\sum_{i=k-r+1}^{k}
\mu_{i}(1)^{2}}{r^{n}\mu_{k}(1)^{n}} \leq
\frac{1}{r^{n-1}\mu_{k}(1)^{n-2}}.$$ Hence
$$|G| < \sum_{i=k-r+1}^{k}
\mu_{i}(1)^{2} + \frac{\sum_{i=k-r+1}^{k}
\mu_{i}(1)^{2}}{r\mu_{k}(1)} + \frac{1}{r}\left( 1 -
\frac{1}{r\mu_{k}(1)}\right)^{-1}$$
$$ = \sum_{i=k-r+1}^{k}
\mu_{i}(1)^{2} + \frac{\sum_{i=k-r+1}^{k}
\mu_{i}(1)^{2}}{r\mu_{k}(1)} + \frac{1}{r - \frac{1}{\mu_{k}(1)}}.$$

\medskip
Now $r > 1$ and also $\mu_{k}(1)
> 1$ (as $G$ is non-Abelian).

\medskip
Thus $$|G| <  1 + \mu_{k}(1) + \sum_{i=k-r+1}^{k} \mu_{i}(1)^{2},$$
so that we now have $$|G| = \mu_{k}(1) + \sum_{i=k-r+1}^{k}
\mu_{i}(1)^{2},$$ as we already have $$|G| \geq \mu_{k}(1) +
\sum_{i=k-r+1}^{k} \mu_{i}(1)^{2}.$$

\medskip
We also now have $$ \frac{|G|}{\mu_{k}(1)} - 1 \leq \theta(1) \leq
\sigma(G^{\#},\theta) \leq \frac{|G|}{\mu_{k}(1)} - 1,$$ so that
$$\theta(1) = \sigma(G^{\#},\theta) = \frac{|G|}{\mu_{k}(1)} - 1.$$
Furthermore, we have $\mu_{i}(1) = \mu_{k}(1)$ whenever $i > k-r.$
Hence $$|G| = \mu_{k}(1) + r\mu_{k}(1)^{2}.$$

\medskip
Furthermore, we also have $|\theta(x)|^{2} = 1 = -\theta(x)$ for all
$x \in G^{\#}$ such that $\theta(x) \neq 0.$ Thus the claim of the
first paragraph of the proof is complete, and in the discussion
which follows, we allow the possibility that $r = 1.$

\medskip
 The orthogonality relations (and the fact that $\mu_{k}(1) =
\sum_{i=1}^{k-r} \mu_{i}(1)^{2}$) now yield that $x \in
\bigcap_{i=1}^{k-r} {\rm ker} \mu_{i}$ whenever $\theta(x) \neq 0.$
Conversely, if $\theta(x) = 0,$ we have $\sum_{i=1}^{k-r}
\mu_{i}(1)\mu_{i}(x) = 0.$

\medskip
Set $K = \bigcap_{i=1}^{k-r} {\rm ker} \mu_{i}.$ Now $$|G| =
\mu_{k}(1) + r\mu_{k}(1)^{2}$$ and $$|K| = 1 + \theta(1) = 1 +
r\mu_{k}(1).$$ Thus $K$ is a Hall subgroup of $G,$ and there is a
complement, $H$ say, to $K$ in $G.$ Since $\mu_{k}(1) = |H|,$ we see
that, for all$i > k-r,$ we have $\mu_{i}(y) = 0$ whenever ${\rm
gcd}(|\langle y \rangle|,|H|) > 1.$ In particular (using a
well-known observation of Feit and Thompson),
$$|C_{G}(y)| = |C_{G/K}(yK)| = |C_{H}(y)|$$ for all $y \in H^{\#}.$
Thus $C_{G}(y) \leq H$ for all $y \in H^{\#}$ and $G$ is a Frobenius
group with kernel $K,$ complement $H.$  Furthermore, every
irreducible character of $G$ which does not contain $K$ in its
kernel has degree $\mu_{k}(1) = [G:K],$ so that every non-trivial
irreducible character of $K$ is linear. Hence $K$ is Abelian.

\medskip
\noindent {\bf PROOF OF COROLLARY 2:} As noted in [1], if $\theta =
\sum_{i=1}^{k} a_{i}\mu_{i},$ then $(\theta,\theta) = \sum_{i<j}
(a_{i}\mu_{j}(1) - a_{j}\mu_{i}(1))^{2}.$

\medskip
\noindent{\bf PROOF OF COROLLARY 3:} By the orthogonality relations
for $G$ and Theorem 1 ( applied to $H$), and noting the presence of
the trivial character, we have
$$\sum_{h \in H^{\#}} |C_{G}(h)| \geq (|H|-1) + (k_{+}(G,H) -
1)\left(\frac{|H|}{\mu(1)} - 1 \right).$$

\medskip
\begin{center}
{\bf REFERENCES}
\end{center}

\medskip
\noindent [1] Robinson,G.R., \emph{ A bound on norms of generalized
characters with applications}, Journal of Algebra, {\bf
212},(1999),660-668.

\medskip
\noindent [2] Robinson,G.R., \emph{ More bounds on norms of
generalized characters with applications to $p$-local bounds and
blocks}, Bulletin of London Mathematical Society, {\bf
37},(4),(2005),555-565.

\medskip
\noindent [3] Robinson,G.R., \emph{ On generalized characters of
nilpotent groups}, Journal of Algebra,{\bf 308},(20007), 822-827.
\end{document}